\documentclass[a4paper,12pt]{article}
\usepackage{amssymb}

\def\smallDeterminant(#1,#2,#3,#4){\left\vert\matrix{#1&#2\cr#3&#4\cr}\right\vert}

\newtheorem{theorem}{Theorem}

\begin{document}

\author{Thabet Abdeljawad\\
\small{\c{C}ankaya University,}\\
\small{Department of Mathematics {\&} Computer Science,}\\
\small{\"{O}gretmenler Cad. 14 06530, Balgat -- Ankara, Turkey}\\
\small{e-mail address: thabet@cankaya.edu.tr}\\
Dumitru B\u{a}leanu\\
\small{\c{C}ankaya University,}\\
\small{Department of Mathematics {\&} Computer Science,}\\
\small{\"{O}gretmenler Cad. 14 06530, Balgat -- Ankara, Turkey}\\
\small{e-mail address: dumitru@cankaya.edu.tr}\\
Fahd Jarad\\
\small{\c{C}ankaya University}\\
\small{Department of Mathematics {\&} Computer Science,}\\
\small{\"{O}gretmenler Cad. 14 06530, Balgat -- Ankara, Turkey}\\
\small{e-mail address: fahd@cankaya.edu.tr}\\
Octavian G. Mustafa\\
\small{\c{C}ankaya University,}\\
\small{Department of Mathematics {\&} Computer Science,}\\
\small{\"{O}gretmenler Cad. 14 06530, Balgat -- Ankara, Turkey}\\
\small{e-mail address: octawian@yahoo.com}\\
Juan J. Trujillo\\
\small{University of La Laguna}\\
\small{Department of Mathematical Analysis,}\\
\small{38271 La Laguna-Tenerife, Spain}\\
\small{e-mail address: jtrujill@ullmat.es}
}

\title{A Fite type result for sequential fractional differential equations}
\date{}
\maketitle

\noindent{\bf Abstract} Given the solution $f$ of the sequential fractional differential equation $_{a}D_{t}^{\alpha}(_{a}D_{t}^{\alpha}f)+P(t)f=0$, $t\in[b,c]$, where $-\infty<a<b<c<+\infty$, $\alpha\in\left(\frac{1}{2},1\right)$ and $P:[a,+\infty)\rightarrow[0,P_{\infty}]$, $P_{\infty}<+\infty$, is continuous, assume that there exist $t_1,t_2\in[b,c]$ such that $f(t_1)=(_{a}D_{t}^{\alpha}f)(t_2)=0$. Then, we establish here a positive lower bound for $c-a$ which depends solely on $\alpha,P_{\infty}$. Such a result might be useful in discussing disconjugate fractional differential equations and fractional interpolation, similarly to the case of (integer order) ordinary differential equations.


\section{Introduction}

The fractional calculus and its plethora of techniques for the study of fractional mathematical models in natural and social disciplines constitute already a vast part of modern science. See the presentations from the authoritative monographs \cite{SamkoKilbas,kilbas_book,Sabatier,MillerRoss,podlubny}. The oscillation/disconjugacy theory of the fractional differential equations, however, is still at the very beginning, cf. \cite{MillerRoss,podlubny,kilbas_book}. To give an example, there are no fractional counterparts of the fundamental comparison and separation theorems of Sturm (\cite{sturm_c1,reid_c1}). 

Complications such as the Fa\`{a} di Bruno formula for the fractional derivative of composite functions or the Leibniz series of the fractional derivative of a product \cite[pp. 96, 98]{podlubny} make such enterprises extremely difficult.

In 1918, the following striking result of W.B. Fite \cite[Theorem VII]{Fite} has opened the way for the disconjugacy analysis \cite{coppel_c1} of (integer order) linear ordinary differential equations: given the real numbers $a<b<c$ and the differential equation
\begin{eqnarray*}
x^{\prime\prime}+P(t)x=0,\quad t\in[b,c],
\end{eqnarray*}
where $P:[a,+\infty)\rightarrow[0,P_{\infty}]$, $P_{\infty}<+\infty$, is a continuous coefficient, assume that $x(t_1)=x^{\prime}(t_2)=0$ for some (non-trivial) solution $x$. Then,
\begin{eqnarray*}
(c-b)\cdot\max\{1,P_{\infty}\}\geq1.
\end{eqnarray*}

In this note, we shall establish a variant of Fite's theorem for the case of linear sequential fractional (ordinary) differential equations (SFDE's). In our result, the estimate of $c-b$ is replaced with an estimate of $c-a$.

To set the hypotheses, fix $\zeta\in(0,1)$. Introduce the Riemann-Liouville (fractional) derivative \cite[p. 68]{podlubny}
\begin{eqnarray*}
(_{a}D_{t}^{\zeta}f)(t)=\frac{1}{\Gamma(1-\zeta)}\cdot\frac{d}{dt}\left[\int_{a}^{t}\frac{f(s)}{(t-s)^{\zeta}}ds\right],\quad t>a,
\end{eqnarray*}
and the Riemann-Liouville (fractional) integral \cite[p. 65]{podlubny} (sometimes denoted with $_{a}I_{t}^{1-\zeta}$, see \cite{kilbas_book})
\begin{eqnarray*}
(_{a}D_{t}^{\zeta-1}f)(t)=\frac{1}{\Gamma(1-\zeta)}\cdot\int_{a}^{t}\frac{f(s)}{(t-s)^{\zeta}}ds,\quad t>a.
\end{eqnarray*}
Obviously, $_{a}D_{t}^{0}f=f$ (a limit case) and $_{a}D_{t}^{\zeta}f=\frac{d}{dt}\left(_{a}D_{t}^{\zeta-1}f\right)$. We shall discuss in the next section the regularity assumptions regarding $f$ that are needed in our investigation.

Consider the following linear sequential differential equation (SFDE)
\begin{eqnarray}
_{a}D_{t}^{\alpha}(_{a}D_{t}^{\alpha}f)+P(t)f=0,\quad t\in[b,c],\label{main_frac}
\end{eqnarray}
where $\alpha\in\left(\frac{1}{2},1\right)$ and the coefficient $P:[a,+\infty)\rightarrow[0,P_{\infty}]$, $P_{\infty}<+\infty$, is assumed continuous. The equation (\ref{main_frac}) is equivalent to the next linear SFDE system
\begin{eqnarray}
\left\{
\begin{array}{ll}
_{a}D_{t}^{\alpha}f=g,\\
_{a}D_{t}^{\alpha}g=-P(t)f,
\end{array}
\right.\quad t\in[b,c],\label{part_syst}
\end{eqnarray}
which is a particular case of the linear SFDE system
\begin{eqnarray}
\left\{
\begin{array}{ll}
(_{a}D_{t}^{\alpha}f)(t)=G(t)g(t)+Q(t),\\
(_{a}D_{t}^{\alpha}g)(t)=R(t)f(t)+V(t),
\end{array}
\right.\quad t\in[b,c],\label{main_sys}
\end{eqnarray}
where the coefficients $G,Q,R,V:[a,+\infty)\rightarrow\mathbb{R}$ are continuous and $G,R$ are bounded: $\Vert G\Vert_{\infty}+\Vert R\Vert_{\infty}<+\infty$.

Assume that, given the solutions $(f_i,g_i)$, with $i=1,2$, of the system (\ref{main_sys}), there exist the points $t_i\in[b,c]$ such that
\begin{eqnarray}
f_1(t_1)=f_2(t_1),\quad g_1(t_2)=g_2(t_2).\label{turi}
\end{eqnarray}
Then, there will exist a positive lower bound for $c-a$ --- this is our main contribution in the present note. It is clear that, in the particular case of (\ref{part_syst}), by taking
\begin{eqnarray}
(f_1,g_1)=(f,{_{a}D_{t}^{\alpha}f})\quad \mbox{and}\quad (f_2,g_2)=(0,0),\label{discoj}
\end{eqnarray}
where $f$ is a non-trivial solution of (\ref{main_frac}), we obtain the fractional version of Fite's theorem.

\section{Preliminaries}

To make the presentation self-contained, we shall discuss here the regularity of the solutions $(f,g)$ and the transformation of (\ref{main_sys}) into an integral system.

\subsection{Regularity of solutions}

Given $\gamma\in(0,1)$, consider the next integral quantity
\begin{eqnarray*}
(Q_{\gamma,A}f)(t)=\int_{a}^{t}\frac{A(s)f(s)}{(t-s)^{\gamma}}ds,\quad t>a,
\end{eqnarray*}
where $A:[a,+\infty)\rightarrow\mathbb{R},f:(a,+\infty)\rightarrow\mathbb{R}$ are continuous and $f$ may be infinite in $t=a$, that is
\begin{eqnarray}
\lim\limits_{t\searrow a}(t-a)^{\gamma}f(t)=f_{a}\in\mathbb{R}.\label{cond_fa}
\end{eqnarray}

In the following, for such functions $f$, we shall use the notations
\begin{eqnarray*}
\Vert f\Vert_{L^{\infty}(a,\gamma;b,c)}=\max\limits_{t\in[b,c]}[(t-a)^{\gamma}\vert f(t)\vert]
\end{eqnarray*}
and
\begin{eqnarray*}
\Vert f\Vert_{L^{\infty}(a,\gamma;c)}=\sup\limits_{t\in(a,c]}[(t-a)^{\gamma}\vert f(t)\vert]
\end{eqnarray*}
for all $c>b>a$. 

Obviously, $\Vert f\Vert_{L^{\infty}(a,\gamma;b,c)}\leq\Vert f\Vert_{L^{\infty}(a,\gamma;c)}$. It is easy to see that
\begin{eqnarray}
\vert f_a\vert\leq\sup\limits_{b\in(a,c)}\Vert f\Vert_{L^{\infty}(a,\gamma;b,c)}=\Vert f\Vert_{L^{\infty}(a,\gamma;c)}.\label{leg_norm}
\end{eqnarray}

Notice that --- for $a<t\leq c$ and $\beta\in(0,1)$ ---
\begin{eqnarray}
\int_{a}^{t}\frac{\vert A(s)f(s)\vert}{(t-s)^{\beta}}ds&=&\int_{a}^{t}\frac{\vert A(s)\vert\cdot (s-a)^{\gamma}\vert f(s)\vert}{(t-s)^{\beta}(s-a)^{\gamma}}ds\nonumber\\
&\leq&\int_{a}^{t}\frac{\Vert A\Vert_{L^{\infty}(a,c)}\Vert f\Vert_{L^{\infty}(a,\gamma;c)}}{(t-s)^{\beta}(s-a)^{\gamma}}ds\nonumber\\
&=&(t-a)^{1-\beta-\gamma}B(1-\gamma,1-\beta)\nonumber\\
&\times&\Vert A\Vert_{L^{\infty}(a,c)}\Vert f\Vert_{L^{\infty}(a,\gamma;c)}.\label{tare_1}
\end{eqnarray}
We have used the change of variables $s=a+\lambda(t-a)$ for $\lambda\in[0,1]$ and the {\it Beta} function $B$, see \cite[p. 6]{podlubny}. Remark also that
\begin{eqnarray*}
\lim\limits_{t\searrow a}(t-a)^{\gamma}(Q_{\beta,A}f)(t)=\lim\limits_{t\searrow a}(t-a)^{\beta}(Q_{\beta,A}f)(t)=0. 
\end{eqnarray*}
In particular, we have the estimate
\begin{eqnarray*}
\vert (Q_{1-\gamma,A}f)(t)\vert&\leq&\int_{a}^{t}\frac{\vert A(s)f(s)\vert}{(t-s)^{1-\gamma}}ds\nonumber\\
&\leq&B(1-\gamma,\gamma)\cdot\Vert A\Vert_{L^{\infty}(a,c)}\Vert f\Vert_{L^{\infty}(a,\gamma;c)}.\label{estim_Q}
\end{eqnarray*}

Let us denote with $C_{\gamma}((a,c],\mathbb{R})$ the set of functions $f\in C((a,c],\mathbb{R})$ which satisfy (\ref{cond_fa}) endowed with the usual function operations. Then, it is easy to check that $X_{\gamma}=(C_{\gamma}((a,c],\mathbb{R}),\Vert\star\Vert_{L^{\infty}(a,\gamma;c)})$ is a Banach space.

So, given $f\in X_{\gamma}$ and $\beta\in(0,1)$, the function $s\longmapsto\frac{A(s)f(s)}{(t-s)^{\beta}}$ belongs to $L^{1}((a,t),\mathbb{R})$ for every $t\in(a,c]$.

The next estimates are crucial for establishing the Fite-type result in this note.

Assume that $a<b\leq t_1<t_2\leq c$. 

We have --- via the change of variables $s=a+\lambda(t_2-a)$ for $\lambda\in[0,1]$ ---
\begin{eqnarray}
&&\int_{t_1}^{t_2}\frac{ds}{(t_2-s)^{\beta}(s-a)^{\gamma}}=(t_2-a)^{1-\beta-\gamma}\cdot\int_{\frac{t_1-a}{t_2-a}}^{1}\frac{d\lambda}{(1-\lambda)^{\beta}\lambda^{\gamma}}\nonumber\\
&&=(t_2-a)^{1-\beta-\gamma}\cdot\int_{0}^{1}\frac{1}{(1-\lambda)^{\beta}\lambda^{\gamma}}\cdot\chi_{\left[\frac{t_1-a}{t_2-a},1\right]}(\lambda)d\lambda\nonumber\\
&&=(t_2-a)^{1-\beta-\gamma}\nonumber\\
&&\times\left[\int_{0}^{\frac{1}{2}}\frac{1}{(1-\lambda)^{\beta}\lambda^{\gamma}}\chi_{\left[\frac{t_1-a}{t_2-a},1\right]}(\lambda)d\lambda+\int_{\frac{1}{2}}^{1}\frac{1}{(1-\lambda)^{\beta}\lambda^{\gamma}}\chi_{\left[\frac{t_1-a}{t_2-a},1\right]}(\lambda)d\lambda\right]\nonumber\\
&&=(t_2-a)^{1-\beta-\gamma}\cdot(I_1+I_2),\label{I1I2}
\end{eqnarray}
where $\chi_S$ denotes the characteristic function of a Lebesgue integrable subset $S$ of $\mathbb{R}$. For simplicity, we shall write $\chi$ instead of $\chi_S$ from now on.

Introduce $p,q,v,w>1$ such that $\frac{1}{p}+\frac{1}{q}=\frac{1}{v}+\frac{1}{w}=1$ and $\gamma p,\beta v<1$. 

In the following, we shall estimate the integrals $I_1,I_2$ by means of H\"{o}lder inequality:
\begin{eqnarray}
I_1&\leq&\int_{0}^{\frac{1}{2}}\frac{1}{\left(\frac{1}{2}\right)^{\beta}\lambda^{\gamma}}\chi(\lambda)d\lambda=2^{\beta}\int_{0}^{\frac{1}{2}}\left(\frac{1}{\lambda}\right)^{\gamma}\chi(\lambda)d\lambda\nonumber\\
&\leq&2^{\beta}\left[\int_{0}^{\frac{1}{2}}\left(\frac{1}{\lambda}\right)^{\gamma p}d\lambda\right]^{\frac{1}{p}}\cdot\left\{\int_{0}^{\frac{1}{2}}[\chi(\lambda)]^{q}d\lambda\right\}^{\frac{1}{q}}\nonumber\\
&=&2^{\beta}\cdot\frac{1}{(1-\gamma p)^{\frac{1}{p}}\cdot 2^{\frac{1}{p}-\gamma}}\cdot\left\{\int_{0}^{\frac{1}{2}}[\chi(\lambda)]d\lambda\right\}^{\frac{1}{q}}\nonumber\\
&=&\frac{2^{\beta+\gamma-\frac{1}{p}}}{(1-\gamma p)^{\frac{1}{p}}}\cdot\left\{\int_{0}^{\frac{1}{2}}[\chi(\lambda)]d\lambda\right\}^{\frac{1}{q}}\leq\frac{2^{\beta+\gamma-\frac{1}{p}}}{(1-\gamma p)^{\frac{1}{p}}}\cdot\left\{\int_{0}^{1}[\chi(\lambda)]d\lambda\right\}^{\frac{1}{q}}\nonumber\\
&=&\frac{2^{\beta+\gamma-\frac{1}{p}}}{(1-\gamma p)^{\frac{1}{p}}}\cdot\left(1-\frac{t_1-a}{t_2-a}\right)^{\frac{1}{q}}\leq c(p,\beta,\gamma)\cdot\left(\frac{t_2-t_1}{t_2-a}\right)^{\frac{1}{q}}\label{estim_I1}
\end{eqnarray} 
and
\begin{eqnarray}
I_2&\leq&2^{\gamma}\int_{\frac{1}{2}}^{1}\left(\frac{1}{1-\lambda}\right)^{\beta}\chi(\lambda)d\lambda\nonumber\\
&\leq&2^{\gamma}\left[\int_{\frac{1}{2}}^{1}\frac{d\lambda}{(1-\lambda)^{\beta v}}\right]^{\frac{1}{v}}\cdot\left\{\int_{0}^{1}[\chi(\lambda)]^{w}d\lambda\right\}^{\frac{1}{w}}\nonumber\\
&=&\frac{2^{\gamma}}{(1-\beta v)^{\frac{1}{v}}\cdot 2^{\frac{1}{v}-\beta}}\cdot\left(\frac{t_2-t_1}{t_2-a}\right)^{\frac{1}{w}}\nonumber\\
&=&c(v,\gamma,\beta)\cdot\left(\frac{t_2-t_1}{t_2-a}\right)^{\frac{1}{w}}.\label{estim_I2}
\end{eqnarray}

By taking into account (\ref{estim_I1}), (\ref{estim_I2}) and (\ref{I1I2}), we deduce that
\begin{eqnarray}
&&\int_{t_1}^{t_2}\frac{ds}{(t_2-s)^{\beta}(s-a)^{\gamma}}\leq(t_2-a)^{1-\beta-\gamma}[c(p,\beta,\gamma)+c(v,\gamma,\beta)]\nonumber\\
&&\times\max\left\{\left(\frac{t_2-t_1}{t_2-a}\right)^{\frac{1}{q}},\left(\frac{t_2-t_1}{t_2-a}\right)^{\frac{1}{w}}\right\}\nonumber\\
&&=(t_2-a)^{1-\beta-\gamma}[c(p,\beta,\gamma)+c(v,\gamma,\beta)]\cdot\left(\frac{t_2-t_1}{t_2-a}\right)^{\frac{1}{\max\{q,w\}}}.\label{estim_fundam_prer_1}
\end{eqnarray}

Next, we would like to estimate the quantity
\begin{eqnarray*}
&&\int_{a}^{t_1}\frac{1}{(s-a)^{\gamma}}\cdot\left[\frac{1}{(t_1-s)^{\beta}}-\frac{1}{(t_2-s)^{\beta}}\right]ds,
\end{eqnarray*}
where $\beta,\gamma\in(0,1)$ and $\beta+\gamma\leq1$.

Noticing that
\begin{eqnarray*}
&&B(1-\gamma,1-\beta)\\
&&=\frac{1}{(t_1-a)^{1-\beta-\gamma}}\cdot\int_{a}^{t_1}\frac{ds}{(s-a)^{\gamma}(t_1-s)^{\beta}}\\
&&=\frac{1}{(t_2-a)^{1-\beta-\gamma}}\cdot\int_{a}^{t_2}\frac{ds}{(s-a)^{\gamma}(t_2-s)^{\beta}},
\end{eqnarray*}
we get that
\begin{eqnarray*}
&&\int_{a}^{t_1}\frac{1}{(s-a)^{\gamma}}\cdot\left[\frac{1}{(t_1-s)^{\beta}}-\frac{1}{(t_2-s)^{\beta}}\right]ds\\
&&=\int_{a}^{t_1}\frac{ds}{(s-a)^{\gamma}(t_1-s)^{\beta}}-\int_{a}^{t_2}\frac{ds}{(s-a)^{\gamma}(t_2-s)^{\beta}}\\
&&+\int_{t_1}^{t_2}\frac{ds}{(t_2-s)^{\beta}(s-a)^{\gamma}}=I_3+I_4.
\end{eqnarray*}
The integral $I_4$ was already evaluated, see (\ref{estim_fundam_prer_1}).

We have --- recall that $1-\beta-\gamma\geq0$ ---
\begin{eqnarray*}
I_3=B(1-\gamma,1-\beta)\cdot[(t_1-a)^{1-\beta-\gamma}-(t_2-a)^{1-\beta-\gamma}]\leq0.
\end{eqnarray*}
 So,
\begin{eqnarray}
&&\int_{a}^{t_1}\frac{1}{(s-a)^{\gamma}}\cdot\left[\frac{1}{(t_1-s)^{1-\gamma}}-\frac{1}{(t_2-s)^{1-\gamma}}\right]ds\leq I_4\nonumber\\
&&\leq(t_2-a)^{1-\beta-\gamma}[c(p,\beta,\gamma)+c(v,\gamma,\beta)]\cdot\left(\frac{t_2-t_1}{t_2-a}\right)^{\frac{1}{\max\{q,w\}}}.\label{estim_I3}
\end{eqnarray}

Finally, take $f\in X_{\gamma}$ and $\gamma,\beta\in(0,1)$ such that $\beta+\gamma\leq1$, and observe that
\begin{eqnarray}
&&\vert (Q_{\beta,A}f)(t_1)-(Q_{\beta,A}f)(t_2)\vert\nonumber\\
&&=\left\vert\int_{a}^{t_1}[A(s)\cdot(s-a)^{\gamma}f(s)]\cdot\frac{1}{(s-a)^{\gamma}}\left[\frac{1}{(t_1-s)^{\beta}}-\frac{1}{(t_2-s)^{\beta}}\right]ds\right.\nonumber\\
&&-\left.\int_{t_1}^{t_2}[A(s)\cdot(s-a)^{\gamma}f(s)]\cdot\frac{ds}{(t_2-s)^{\beta}(s-a)^{\gamma}}\right\vert\nonumber\\
&&\leq\Vert A\Vert_{L^{\infty}(a,c)}\Vert f\Vert_{L^{\infty}(a,\gamma;b,c)}\cdot\int_{t_1}^{t_2}\frac{ds}{(t_2-s)^{\beta}(s-a)^{\gamma}}\nonumber\\
&&+\Vert A\Vert_{L^{\infty}(a,c)}\Vert f\Vert_{L^{\infty}(a,\gamma;c)}\int_{a}^{t_1}\frac{1}{(s-a)^{\gamma}}\cdot\left[\frac{1}{(t_1-s)^{\beta}}-\frac{1}{(t_2-s)^{\beta}}\right]ds\nonumber\\
&&\leq\Vert A\Vert_{L^{\infty}(a,c)}\Vert f\Vert_{L^{\infty}(a,\gamma;c)}\cdot2(t_2-a)^{1-\beta-\gamma}[c(p,\beta,\gamma)+c(v,\gamma,\beta)]\nonumber\\
&&\times\left(\frac{t_2-t_1}{t_2-a}\right)^{\frac{1}{\max\{q,w\}}}\nonumber\\
&&=2(t_2-a)^{1-\beta-\gamma-\frac{1}{\max\{q,w\}}}\cdot[c(p,\beta,\gamma)+c(v,\gamma,\beta)]\nonumber\\
&&\times(t_2-t_1)^{\frac{1}{\max\{q,w\}}}\cdot\Vert A\Vert_{L^{\infty}(a,c)}\Vert f\Vert_{L^{\infty}(a,\gamma;c)}\nonumber\\
&&=C\cdot(t_2-a)^{1-\beta-\gamma-\frac{1}{\max\{q,w\}}}\nonumber\\
&&\times(t_2-t_1)^{\frac{1}{\max\{q,w\}}}\cdot\Vert A\Vert_{L^{\infty}(a,c)}\Vert f\Vert_{L^{\infty}(a,\gamma;c)},\label{estim_cont}
\end{eqnarray}
where $C=C(p,v,\beta,\gamma)=2[c(p,\beta,\gamma)+c(v,\gamma,\beta)]$.

Before going further, recall the following inequality: given the numbers $A,B\geq0$ and $\beta\in(0,1)$, one has
\begin{eqnarray}
(A+B)^{\beta}\leq A^{\beta}+B^{\beta}.\label{ineq_abbeta}
\end{eqnarray}
For $A=t_1-a$ and $B=t_2-t_1$, the inequality (\ref{ineq_abbeta}) leads to
\begin{eqnarray*}
0\leq(t_{2}-a)^{\beta}-(t_{1}-a)^{\beta}\leq(t_{2}-t_{1})^{\beta},
\end{eqnarray*}
an estimate that will be used in the next series of computations.

Now,
\begin{eqnarray}
&&\vert(t_1-a)^{\beta}(Q_{\beta,A}f)(t_1)-(t_2-a)^{\beta}(Q_{\beta,A}f)(t_2)\vert\nonumber\\
&&=\vert(t_1-a)^{\beta}[(Q_{\beta,A}f)(t_1)-(Q_{\beta,A}f)(t_2)]\nonumber\\
&&-(Q_{\beta,A}f)(t_2)[(t_2-a)^{\beta}-(t_1-a)^{\beta}]\vert\nonumber\\
&&\leq(t_1-a)^{\beta}\cdot\vert (Q_{\beta,A}f)(t_1)-(Q_{\beta,A}f)(t_2)\vert\nonumber\\
&&+\vert(Q_{\beta,A}f)(t_2)\vert\cdot(t_2-t_1)^{\beta}\nonumber\\
&&\leq(t_1-a)^{\beta}\cdot C(t_2-a)^{1-\beta-\gamma-\frac{1}{\max\{q,w\}}}(t_2-t_1)^{\frac{1}{\max\{q,w\}}}\nonumber\\
&&\times\Vert A\Vert_{L^{\infty}(a,c)}\Vert f\Vert_{L^{\infty}(a,\gamma;c)}\nonumber\\
&&+(t_2-a)^{1-\beta-\gamma}B(1-\gamma,1-\beta)\cdot\Vert A\Vert_{L^{\infty}(a,c)}\Vert f\Vert_{L^{\infty}(a,\gamma;c)}\nonumber\\
&&\times(t_2-t_1)^{\beta}\nonumber\\
&&\leq\left[C(c-a)^{\beta}\cdot(t_2-a)^{1-\beta-\gamma-\frac{1}{\max\{q,w\}}}\right.\label{estim_spec_00}\\
&&\left.+(c-a)^{1-\beta-\gamma}\cdot B(1-\gamma,1-\beta)\right]\Vert A\Vert_{L^{\infty}(a,c)}\Vert f\Vert_{L^{\infty}(a,\gamma;c)}\nonumber\\
&&\times\max\left\{(t_2-t_1)^{\frac{1}{\max\{q,w\}}},(t_2-t_1)^{\beta}\right\}\nonumber\\
&&\leq D\Vert A\Vert_{L^{\infty}(a,c)}\Vert f\Vert_{L^{\infty}(a,\gamma;c)}\max\left\{(t_2-t_1)^{\frac{1}{\max\{q,w\}}},(t_2-t_1)^{\beta}\right\},\label{estim_cont_2}
\end{eqnarray}
where $D=D(p,v,\beta,\gamma,b-a,c-a)$. We took into account the estimates (\ref{tare_1}), (\ref{estim_cont}).

Let us make a comment regarding $D$ and (\ref{estim_spec_00}). If we don't know the sign of the quantity $1-\beta-\gamma-\frac{1}{\max\{q,w\}}$ then we are forced to use the raw estimate
\begin{eqnarray*}
(t_2-a)^{1-\beta-\gamma-\frac{1}{\max\{q,w\}}}\leq(b-a)^{1-\beta-\gamma-\frac{1}{\max\{q,w\}}}+(c-a)^{1-\beta-\gamma-\frac{1}{\max\{q,w\}}}
\end{eqnarray*}
which shall reflect on the size of $D$. If, on the other hand, we have $1-\beta-\gamma-\frac{1}{\max\{q,w\}}>0$ then
\begin{eqnarray*}
(t_2-a)^{1-\beta-\gamma-\frac{1}{\max\{q,w\}}}\leq(c-a)^{1-\beta-\gamma-\frac{1}{\max\{q,w\}}},\quad t_2\in[b,c],
\end{eqnarray*}
and so
\begin{eqnarray}
D=C\cdot(c-a)^{1-\gamma-\frac{1}{\max\{q,w\}}}+(c-a)^{1-\beta-\gamma}\cdot B(1-\gamma,1-\beta).\label{estim_spec_01}
\end{eqnarray}
Observe that $D$ is {\it independent} of $b$, a fact of crucial importance for our investigation.

In conclusion, given $f\in X_{\gamma}$, the functions $t\longmapsto(Q_{\beta,A}f)(t)$ and $t\longmapsto(t-a)^{\beta}(Q_{\beta,A}f)(t)$, with $0<\beta\leq1-\gamma$, are continuous in $(a,c]$ by being uniformly continuous on each interval $[b,c]$, with $b\in(a,c)$.

In this note, by a {\it solution} of the fractional differential system (\ref{main_sys}) we mean any pair $(f,g)$ of functions from $X_{1-\alpha}$ which verify the equations of the system.

\subsection{Integral representation of (\ref{main_sys})}

Consider now the first of equations (\ref{main_sys}), namely
\begin{eqnarray}
(_{a}D_{t}^{\alpha}f)(t)=G(t)g(t)+Q(t),\quad t\in[b,c],\label{gen_eq_1}
\end{eqnarray}
where $\alpha\in\left(\frac{1}{2},1\right)$.

To use the estimates from the preceding subsection, introduce the quantities
\begin{eqnarray*}
0<\gamma=\beta=1-\alpha<\frac{1}{2}.
\end{eqnarray*}
There exist $p=v>1$ such that
\begin{eqnarray*}
\gamma p=\beta v<\frac{1}{2}.
\end{eqnarray*}
Also, $q=w$.

An important consequence of this choice of constants is that
\begin{eqnarray*}
1-\beta-\gamma-\frac{1}{\max\{q,w\}}=1-2\gamma-\frac{1}{q}=\frac{1}{p}-2\gamma>0,
\end{eqnarray*}
leading to --- recall (\ref{estim_spec_01}) --- a quantity $D$ which is independent of $b\in(a,c)$.

Also, since $c(p,\gamma,\gamma)=\frac{2^{2\gamma-\frac{1}{p}}}{(1-\gamma p)^{\frac{1}{p}}}<\frac{2^{2\gamma-\frac{1}{p}}}{\left(\frac{1}{2}\right)^{\frac{1}{p}}}=2^{2\gamma}$, we have
\begin{eqnarray*}
C=2[c(p,\beta,\gamma)+c(v,\gamma,\beta)]=4c(p,\gamma,\gamma)<2^{2(1+\gamma)}=2^{2(2-\alpha)}
\end{eqnarray*}
and respectively
\begin{eqnarray}
D&<&2^{2(2-\alpha)}(c-a)^{\alpha-\frac{1}{q}}+B(\alpha,\alpha)(c-a)^{2\alpha-1}\nonumber\\
&\leq&[2^{2(2-\alpha)}+B(\alpha,\alpha)]\cdot\frac{(c-a)^{\alpha}}{\min\{(c-a)^{\frac{1}{q}},(c-a)^{1-\alpha}\}}.\label{trD_1}
\end{eqnarray}

Now, recasting (\ref{gen_eq_1}) as
\begin{eqnarray*}
\frac{d}{dt}\left[\int_{a}^{t}\frac{f(s)}{(t-s)^{\alpha}}ds\right]=\Gamma(1-\alpha)[G(t)g(t)+Q(t)],
\end{eqnarray*}
we deduce that
\begin{eqnarray}
&&\int_{a}^{t}\frac{f(s)}{(t-s)^{\alpha}}ds-\lim\limits_{t\searrow a}\left[\int_{a}^{t}\frac{f(s)}{(t-s)^{\alpha}}ds\right]\nonumber\\
&&=\Gamma(1-\alpha)\int_{a}^{t}[G(s)g(s)+Q(s)]ds.\label{integr_interm1}
\end{eqnarray}

By taking into account the simple estimate
\begin{eqnarray*}
&&\left\vert\int_{a}^{t}\frac{1}{(t-s)^{\alpha}(s-a)^{1-\alpha}}\cdot[(s-a)^{1-\alpha}f(s)-f_a]ds\right\vert\\
&&\leq B(\alpha,1-\alpha)\cdot\sup\limits_{s\in(a,t]}\vert(s-a)^{1-\alpha}f(s)-f_a\vert,
\end{eqnarray*}
the identity (\ref{integr_interm1}) yields
\begin{eqnarray*}
\int_{a}^{t}\frac{f(s)}{(t-s)^{\alpha}}ds-f_{a}B(\alpha,1-\alpha)=\Gamma(1-\alpha)\int_{a}^{t}[G(s)g(s)+Q(s)]ds.
\end{eqnarray*}

Further, we have
\begin{eqnarray*}
&&\int_{a}^{x}\frac{1}{(x-t)^{1-\alpha}}\int_{a}^{t}\frac{f(s)}{(t-s)^{\alpha}}dsdt-f_{a}B(\alpha,1-\alpha)\cdot\frac{(x-a)^{\alpha}}{\alpha}\\
&&=\Gamma(1-\alpha)\int_{a}^{x}\frac{1}{(x-t)^{1-\alpha}}\int_{a}^{t}[G(s)g(s)+Q(s)]dsdt.
\end{eqnarray*}

Next, the left hand member of the preceding identity becomes
\begin{eqnarray*}
&&\int_{a}^{x}f(s)\int_{s}^{x}\frac{dt}{(x-t)^{1-\alpha}(t-s)^{\alpha}}ds-f_{a}B(\alpha,1-\alpha)\cdot\frac{(x-a)^{\alpha}}{\alpha}\\
&&=B(\alpha,1-\alpha)\cdot\int_{a}^{x}f(s)ds-f_{a}B(\alpha,1-\alpha)\cdot\frac{(x-a)^{\alpha}}{\alpha},
\end{eqnarray*}
while the right hand member reads as, by means of the Abel computations \cite[p. 32, eq. (2.13)]{SamkoKilbas},
\begin{eqnarray*}
\Gamma(1-\alpha)\int_{a}^{x}\int_{a}^{t}\frac{G(s)g(s)+Q(s)}{(t-s)^{1-\alpha}}dsdt.
\end{eqnarray*}

By differentiating with respect to $x$, where $x>a$, the identity
\begin{eqnarray*}
&&B(\alpha,1-\alpha)\cdot\int_{a}^{x}f(s)ds-f_{a}B(\alpha,1-\alpha)\cdot\frac{(x-a)^{\alpha}}{\alpha}\\
&&=\Gamma(1-\alpha)\int_{a}^{x}\int_{a}^{t}\frac{G(s)g(s)+Q(s)}{(t-s)^{1-\alpha}}dsdt,
\end{eqnarray*}
we get
\begin{eqnarray}
f(x)=\frac{f_a}{(x-a)^{1-\alpha}}+\frac{1}{\Gamma(\alpha)}\int_{a}^{x}\frac{G(s)g(s)+Q(s)}{(x-s)^{1-\alpha}}ds,\quad x\in(a,c].\label{interm_2}
\end{eqnarray}

Notice that, since $g\in X_{1-\alpha}$, the function $s\longmapsto G(s)g(s)+Q(s)$, denoted $H$, is also in $X_{1-\alpha}$. We have $H_a=G(a)\cdot g_a$. As $\beta+\gamma=(1-\alpha)+(1-\alpha)<1$, the integral from the right side of (\ref{interm_2}) is $(Q_{\beta,1}H)(x)$ and both $(Q_{\beta,1}H)(x),(x-a)^{\beta}(Q_{\beta,1}H)(x)$ are continuous in $(a,c]$.

Similarly,
\begin{eqnarray}
g(x)=\frac{g_a}{(x-a)^{1-\alpha}}+\frac{1}{\Gamma(\alpha)}\int_{a}^{x}\frac{R(s)f(s)+V(s)}{(x-s)^{1-\alpha}}ds,\quad x\in(a,c].\label{interm_3}
\end{eqnarray}

For other derivations of formulas (\ref{interm_2}), (\ref{interm_3}), see \cite[p. 837 and following]{SamkoKilbas} and \cite[pp. 122--123]{podlubny}

\section{Main result}

\begin{theorem}\label{frac_Fite} Set $p>1$ such that $(1-\alpha)p<\frac{1}{2}$ and take $q=\frac{p-1}{p}$. Assume that there exist the solutions $(f_i,g_i)$ of the SFDE system (\ref{main_sys}) such that (\ref{turi}) holds for some $t_1,t_2\in[b,c]$. Then,
\begin{eqnarray*}
m\cdot(c-a)^{\alpha}\frac{\max\left\{(c-a)^{\frac{1}{q}},(c-a)^{1-\alpha}\right\}}{\min\{(c-a)^{\frac{1}{q}},(c-a)^{1-\alpha}\}}\geq\frac{\Gamma(\alpha)}{2^{2(2-\alpha)}+B(\alpha,\alpha)},
\end{eqnarray*}
where $m=\max\{\Vert G\Vert_{\infty},\Vert R\Vert_{\infty}\}$.
\end{theorem}

\textbf{Proof.} We have the relations
\begin{eqnarray*}
\left\{
\begin{array}{ll}
f_i(t)=\frac{(f_{i})_a}{(t-a)^{1-\alpha}}+\frac{1}{\Gamma(\alpha)}\int_{a}^{t}\frac{G(s)g_i(s)+Q(s)}{(t-s)^{1-\alpha}}ds,\\
g_i(t)=\frac{(g_{i})_a}{(t-a)^{1-\alpha}}+\frac{1}{\Gamma(\alpha)}\int_{a}^{t}\frac{R(s)f_i(s)+V(s)}{(t-s)^{1-\alpha}}ds,
\end{array}
\right.\quad t\in[b,c],\quad i\in\{1,2\}.
\end{eqnarray*}

Take $b\leq t\leq t_1\leq c$. By means of (\ref{estim_cont_2}), we deduce that
\begin{eqnarray*}
&&\vert(t-a)^{1-\alpha}(f_2-f_1)(t)\vert\\
&&=\vert (t_1-a)^{1-\alpha}(f_2-f_1)(t_1)-(t-a)^{1-\alpha}(f_2-f_1)(t)\vert\\
&&=\left\vert\frac{1}{\Gamma(\alpha)}\cdot[(t-a)^{1-\alpha}(Q_{1-\alpha,G}(g_2-g_1))(t)\right.\\
&&-\left.(t_1-a)^{1-\alpha}(Q_{1-\alpha,G}(g_2-g_1))(t_1)]\right\vert\\
&&\leq\frac{D}{\Gamma(\alpha)}\cdot\Vert G\Vert_{L^{\infty}(a,c)}\Vert g_2-g_1\Vert_{L^{\infty}(a,1-\alpha;c)}\\
&&\times\max\left\{(t_1-t)^{\frac{1}{q}},(t_1-t)^{1-\alpha}\right\}\\
&&\leq\frac{D\Vert G\Vert_{\infty}}{\Gamma(\alpha)}\max\left\{(c-b)^{\frac{1}{q}},(c-b)^{1-\alpha}\right\}\cdot\Vert g_2-g_1\Vert_{L^{\infty}(a,1-\alpha;c)}\\
&&\leq\frac{D\Vert G\Vert_{\infty}}{\Gamma(\alpha)}\max\left\{(c-a)^{\frac{1}{q}},(c-a)^{1-\alpha}\right\}\cdot\Vert g_2-g_1\Vert_{L^{\infty}(a,1-\alpha;c)}.
\end{eqnarray*}

The case $b\leq t_1\leq t\leq c$ leads to the same conclusion. So,
\begin{eqnarray}
&&\Vert f_2-f_1\Vert_{L^{\infty}(a,1-\alpha;b,c)}\nonumber\\
&&\leq E\cdot\Vert G\Vert_{\infty}\cdot\Vert g_2-g_1\Vert_{L^{\infty}(a,1-\alpha;c)}\nonumber\\
&&\leq E\cdot\max\{\Vert G\Vert_{\infty},\Vert R\Vert_{\infty}\}\cdot\Vert g_2-g_1\Vert_{L^{\infty}(a,1-\alpha;c)},\label{X_1}
\end{eqnarray}
where, via (\ref{trD_1}),
\begin{eqnarray}
&&E=\frac{D}{\Gamma(\alpha)}\max\left\{(c-a)^{\frac{1}{q}},(c-a)^{1-\alpha}\right\}\nonumber\\
&&\leq\frac{2^{2(2-\alpha)}+B(\alpha,\alpha)}{\Gamma(\alpha)}\cdot(c-a)^{\alpha}\frac{\max\left\{(c-a)^{\frac{1}{q}},(c-a)^{1-\alpha}\right\}}{\min\{(c-a)^{\frac{1}{q}},(c-a)^{1-\alpha}\}}.\label{label_E}
\end{eqnarray}

Similarly,
\begin{eqnarray}
&&\Vert g_2-g_1\Vert_{L^{\infty}(a,1-\alpha;b,c)}\nonumber\\
&&\leq E\cdot\Vert R\Vert_{\infty}\cdot\Vert f_2-f_1\Vert_{L^{\infty}(a,1-\alpha;c)}\nonumber\\
&&\leq E\cdot\max\{\Vert G\Vert_{\infty},\Vert R\Vert_{\infty}\}\cdot\Vert f_2-f_1\Vert_{L^{\infty}(a,1-\alpha;c)}.\label{X_2}
\end{eqnarray}

In conclusion, as $E$ is independent of $b$, by taking into account (\ref{X_1}), (\ref{X_2}) and (\ref{leg_norm}), we get
\begin{eqnarray*}
X\leq E\cdot\max\{\Vert G\Vert_{\infty},\Vert R\Vert_{\infty}\}\cdot X,
\end{eqnarray*}
where $X=\max\{\Vert f_2-f_1\Vert_{L^{\infty}(a,1-\alpha;c)},\Vert g_2-g_1\Vert_{L^{\infty}(a,1-\alpha;c)}\}$. 

The proof is complete. $\square$

\section{Particular cases}

In the particular case of (\ref{main_frac}),(\ref{discoj}), Theorem \ref{frac_Fite} provides an upper bound to the length of the open interval $(a,c)$ in which at least one of $f$ and $_{a}D_{t}^{\alpha}f$ has no zero. Such an interval can be called a (fractional) {\it disconjugacy} interval of SFDE (\ref{main_frac}). This type of investigation might be useful for problems related to the localization of zeros of various fractional variants of sine and cosine type of functions, cf. \cite{luchko}, \cite[p. 19]{podlubny}.

Another particular case is that of sequential relaxation-oscillation equations, namely
\begin{eqnarray*}
_{a}D_{t}^{\alpha}(_{a}D_{t}^{\alpha}f)+P\cdot f=V(t),\quad t\in[b,c],
\end{eqnarray*}
where $\alpha\in\left(\frac{1}{2},1\right)$, the coefficient $P>0$ is constant and the free term $V:[a,+\infty)\rightarrow\mathbb{R}$ is continuous. As before, we are interested in the solutions $f$ such that $f,{\;}_{a}D_{t}^{\alpha}f\in X_{1-\alpha}$. Such solutions exhibit several zeros for various terms $V$ --- see the numerical investigations from \cite[pp. 224--226, 238--242]{podlubny}.

\section{Acknowledgments.}
The work on this paper has been done during the visit of O.G. Mustafa to \c{C}ankaya University in the fall of 2008. He is grateful to the people from the Department of Mathematics and Computer Science for the friendly and enthusiastic working atmosphere. J.J. Trujillo has been supported by the MICINN of Spain via MTM2007-60246.

\end{document}